\documentclass[12pt]{amsart}
\usepackage{amsmath,amssymb,amsthm,upref,graphicx,mathrsfs}

\usepackage{color}
\usepackage[
  colorlinks=true,
  linkcolor=blue,
  citecolor=blue,
  urlcolor=blue]{hyperref}

\numberwithin{equation}{section}

\newtheorem{thm}{Theorem}[section]

\newtheorem{cor}[thm]{Corollary}
\newtheorem{lem}[thm]{Lemma}
\newtheorem{defi}[thm]{Definition}
\newtheorem{remark}[thm]{Remark}
\newtheorem{example}[thm]{Example}
\newtheorem{pb}[thm]{Problem}
\newenvironment{rk}{\begin{remark}\rm}{\end{remark}}

\newcommand{\R}{{\mathbb R}}

\newcommand{\M}{{\mathcal M}}
\newcommand{\E}{{\mathbb E}}

\newcommand{\bs}{\begin{split}}
\newcommand{\es}{\end{split}}

\newcommand{\be}{\begin{eqnarray*}}
\newcommand{\ee}{\end{eqnarray*}}
\newcommand{\beq}{\begin{equation}}
\newcommand{\eeq}{\end{equation}}

\begin{document}

\title[\bf Noncommutative harmonic analysis on semigroups]
  {\bf Noncommutative harmonic analysis on semigroups}

\authors

\author[Y. Jiao]{Yong Jiao}
\address{Yong Jiao \\ School of Mathematical Sciences,
Central South University, 410083 Changsha, China}
\email{jiaoyong@csu.edu.cn}

\author[M. Wang]{Maofa Wang}
\address{Maofa Wang\\ School of Mathematics and Statistics, Wuhan University, 430072 Wuhan, China}
\email{mfwang.math@whu.edu.cn}

\makeatletter
\renewcommand{\@makefntext}[1]{#1}
\makeatother \footnotetext{\noindent
Y.Jiao (jiaoyong@csu.edu.cn) was supported by NSFC(No.11471337), Hunan NSF(14JJ1004); M.Wang (mfwang.math@whu.edu.cn) was supported by NSFC(No.11271293,11431011,11471251).}

\keywords{Semigroup, Noncommutative
$L_p$-space, Transference, Noncommutative Markov dilation.}
\subjclass[2010]{Primary 46L52; Secondary 42B25.}

%
%
\begin{abstract} {In this paper we obtain some noncommutative multiplier theorems and
 maximal inequalities on semigroups. As applications,
we obtain the corresponding individual ergodic theorems. Our main results extend
some classical results {of} Stein and Cowling {on one hand, and simplify the main arguments of Junge-Le Merdy-Xu's related work \cite{JLX}.}}
\end{abstract}

\maketitle

%
\section{Introduction and Preliminaries }

 Let $L$ be a densely defined positive operator on $L_2(\Omega)$, where $\Omega$
is a $\sigma$-finite measure space. Suppose that $\{P_\lambda\}$ is the spectral resolution of $L$:
$$Lf=\int_0^\infty\lambda dP_\lambda f,\quad \quad f\in\text{Dom}(L).$$

\noindent If $m$ is a bounded function on $[0,\infty)$, then by the
spectral theorem, the multiplier operator $m(L)$ defined by
$$m(L)f=\int_0^\infty m(\lambda) dP_\lambda f,\quad \quad f\in L_2(\Omega),$$

\noindent is bounded on $L_2(\Omega).$ Let $(T_t)_{t>0}=(e^{-tL})_{t>0}$ be
the operator semigroup, which we always assume  satisfies the contraction property:
$$\|T_tf\|_p\leq\|f\|_p,\quad \quad f\in L_2(\Omega)\cap
L_p(\Omega),$$

\noindent wherever  $1\leq p\leq\infty.$   Stein
\cite{Stein} developed a Littlewood-Paley theory for such
semigroups, with some additional hypotheses. By use of
transference techniques, Coifman and Weiss
\cite{cw}, Cowling \cite{cowling} presented an alternative and
simpler approach to obtain some multiplier results and maximal
inequalities. Indeed, Cowling \cite{cowling} showed that $m(L)$,
originally defined on $L_2(\Omega)$ via the spectral theorem, is
exactly a bounded operator on $L_p(\Omega)$ for $1<p<\infty$,
whenever $m$ has a bounded analytic extension on some sector $\Sigma_\phi$ with $\phi>\pi|1/p-1/2|$, where
$$\Sigma_\phi=\{{z\in\mathbb{C}: |\arg z|}<\phi\}.$$
\noindent  More precisely, the following estimate
holds: \begin{eqnarray}\|m(L)f\|_p\leq
C_{p,\phi}\|m\|_\infty\|f\|_p,\quad \quad f\in
L_p(\Omega).
\end{eqnarray}
In other words, the generators of the
symmetric contraction semigroups have a $H^\infty$
functional calculus on each sector $\Sigma_\phi$ with
$\phi>\pi|1/p-1/2|$.

Recently,
more attention was turned
to diffusion semigroups on noncommutative space $L_p(\mathcal {M})$
associated to a von Neumann algebra, see for instance
\cite{JLX,jm2,jm1,JX,K,m}. In this paper, we consider the semigroup
$(T_t)_{t>0}$ acting on noncommutative $L_p$-space associated with
$(\mathcal {M},\tau)$, where $\mathcal {M}$ is a von Neumann algebra
with a normal {finite} faithful trace $\tau$. Under reasonable
hypotheses, we obtain some noncommutative multiplier theorems,
the noncommutative version of (1.1), which positively
answers the question raised in \cite[Remark 5.9]{JLX}. Namely, the
generators of the noncommutative diffusion semigroups also have a
$H^\infty$ functional calculus on each sector $\Sigma_\phi$ with
$\phi>\pi|1/p-1/2|$. By means of the multiplier theorems, we
establish some noncommutative maximal inequalities and
 individual ergodic theorems. {It is worth pointing out
that the key point of Cowling's method in \cite{cowling}
is to combine the transference technique and Fendler's dilation
theorem.  A noncommutative version of Fendler's dilation theorem
has been recently achieved by Junge-Ricard-Shlyakhtenko
\cite{jrs} (see also Dabrowski \cite{Dab}). Armed with this result,
we can extend Cowling's method to the noncommutative setting.
In this way, we recover the main results of \cite{JLX} by a very simple method.
This is a major advantage of our method over that of \cite{JLX}.}

Now we introduce some preliminaries which will be used in
the sequel. We shall work on a von Neumann algebra
$\mathcal {M}$ equipped with a normal {finite} faithful trace $\tau.$ For
$1\leq p<\infty,$ let $L_p(\mathcal {M},\tau)$ or simply
$L_p(\mathcal {M})$ be the associated noncommutative $L_p$ space.
Namely, $L_p(\mathcal {M})$ is the completion of $\mathcal {{\mathcal {M}}}$ with
the norm $\|x\| _p=(\tau(|x|^p))^{1/p}$, where $|x|=(x^*x)^{1/2}$ is  {the modulus of $x$.} 
By convention we set $L_\infty(\mathcal {M})=\mathcal
{M}$, equipped with the operator norm. Like the commutative $L_p$-spaces, one has the duality:
$L_p(\mathcal {M})^*=L_q(\mathcal {M})$ via $(x,y)\mapsto\tau(xy)$,
for  $x\in L_p{(\mathcal {M})}, y\in L_q{(\mathcal {M})}$ with $1\leq p<\infty$ and $1/p+1/q=1.$ It is also well known that
$L_p(\mathcal {M})$ has UMD property for $1<p<\infty.$ We refer to
\cite{pxhand} for more information and more historical
references on noncommutative $L_p$-spaces.

We say an operator $T$ on $\M$ is completely positive if $T\otimes
I_n$ is positive on $\M\otimes M_n$ for each $n.$ Here, $M_n$ is the
algebra of $n\times n$ matrices and $I_n$ is the identity operator on
$M_n.$ Now we introduce the standard noncommutative semigroup.
That is, $(T_t)$ is a semigroup
of completely positive maps on a {finite} von Neumann algebra $\M$
satisfying the following conditions:

1) Every $T_t$ is normal on $\M$ such that
$T_t(1)=1;$

2) Every $T_t$ is selfadjoint with respect to the trace $\tau$, i.e.
$\tau\big(T_t(x)y\big)=\tau\big(xT_t(y)\big)$;

3) The family $(T_t)$
is strongly continuous, i.e. $\lim _{t\to 0}T_tx=x$ with respect to the
strong {operator} topology in $\M$ for any $x\in \M.$

Let us note that the
first two conditions imply that $\tau(T_tx)=\tau(x)$ for all $x$, so $T_t$  is
faithful and contractive on $L_1(\M)$. By interpolation {technique}, $T_t$ can be
extend to a contraction on $L_p(\mathcal {M})$ for $1\leq p<\infty$
and satisfies $\lim_{t\to 0} T_tx=x$ in $L_p(\mathcal {M})$ for any $x\in
L_p(\mathcal {M})$. Let us recall that such a semigroup admits an
infinitesimal generator $L$, i.e., $T_t=e^{-tL}.$ We
refer to \cite {JX} for more details.

We say that a standard semigroup $(T_t)$ on a {finite} von Neumann
algebra $\M$ admits a Markov dilation if there exists a
larger {finite} von Neumann algebra $\mathcal {N}$, an increasing
filtration $(\mathcal {N}_{s]})_{s\geq0}$ with conditional
expectation $\mathcal {N}_{s]}=\E_{s]}(\mathcal {N})$ and trace
preserving
*-homomorphisms $\pi_s:\M\longrightarrow\mathcal {N}$
such that $\pi_s(\M)\subset\mathcal
{N}_{s]}$ and
$$\E_{s]}\big(\pi_t(x)\big)=\pi_s(T_{t-s}x),
\quad \quad 0\leq s{<}t<\infty,\quad x\in \M.$$
In \cite{jrs}, the authors proved that every semigroup of
completely positive unital selfadjoint maps on a {finite} von
Neumann algebra admits a Markov dilation. Moreover, the authors in \cite{jrs} extended the  Markov dilation above to all of $\mathbb{R}$ by using the ultraproduct argument. Namely, there exists a
new {finite} von Neumann algebra $\mathcal {N}$, an increasing
filtration $(\mathcal {N}_{s]})_{-\infty<s<\infty}$ with conditional
expectation $\mathcal {N}_{s]}=\E_{s]}(\mathcal {N})$ and trace
preserving
*-homomorphisms $\pi_s:\M\longrightarrow\mathcal {N}$
such that $\pi_s(\M)\subset\mathcal
{N}_{s]}$ and
$$\E_{s]}\big(\pi_t(x)\big)=\pi_s(T_{t-s}x),
\quad \quad -\infty< s{<}t<\infty,\quad x\in \M.$$

Our paper is organized as follows. Section 2 is on noncommutative
multiplier theorems. The noncommutative maximal inequalities are
proved in Section 3.  As applications, in Section 4, we
give some individual ergodic theorems.

In the rest of the paper we use the same
letter $C$  to denote various positive constants which may change at
each occurrence. Variables indicating the dependency of constants
$C$ will be often specified in the parenthesis. We use the notation
$X\lesssim Y$  or $Y\gtrsim X$ for nonnegative quantities $X$ and
$Y$ to mean $X\le CY$ for some inessential constant $C>0$.
Similarly, we use the notation $X\approx Y$ if both $X\lesssim Y$
and $Y \lesssim X$ hold.

\section{Noncommutative multiplier theorems}\label{section B}

In this section, {we first give a noncommutative Fourier multiplier
theorem by applying the following Junge-Ricard-Shlyakhtenko dilation theorem \cite[Theorem 5 and Corollary 4.5]{jrs}, which plays a crucial role in our proof.}

{\begin{thm}\label{JXS}  Let $(T_t)$ be a semigroup of completely positive unital and selfadjoint maps on a
finite von Neumann algebra $(\M,\tau ).$ Then $(T_t)$ admits a Markov dilation.
\end{thm}}

{In the sequel, we suppose that the spectral projection $P_0$ onto
the kernel of $L$ is trivial on $L_p(\mathcal
{M})$ and we hence do not need to consider the definition of $m(0)$.}

\begin{thm}
\label{main}  Suppose that $m$ is a bounded holomorphic
function on the sector $\Sigma_{\pi/2}$. Let $\Phi$ be the
distribution on $\mathbb{R}$ whose Fourier transform is the bounded
function defined (almost everywhere) by the
formula$$\hat{\Phi}(v)=m(iv),\quad v\in \mathbb{R},$$ where
$m(iv)$ is the non-tangential limit. If for some $p\in [1,\infty]$
and  all $f$ in $L_p\big(\mathbb{R},L_p(\mathcal {M})\big)$
$$\|\Phi*f\|_{L_p\big(\mathbb{R},L_p(\mathcal {M})\big)}\leq C\|f\|_{L_p\big(\mathbb{R},L_p(\mathcal {M})\big)},$$
for some positive constant $C$, then for all $x\in \mathcal {M},$
\begin{eqnarray}
\|m(L)x\|_p\leq C\|x\|_p.
\end{eqnarray}

\noindent Consequently, $m(L)$ extends uniquely
to a bounded operator on $L_p(\mathcal {M})$, still denoted \textcolor{red}{by} $m(L)$, of norm at most $C$.
\end{thm}

\noindent {\bf Proof} \  We proceed the proof by a standard transference argument.
Since the distribution $\Phi$ is defined as follows:
$$\hat{\Phi}(v)=\int_{-\infty}^{+\infty}e^{-iuv}\Phi(u)du=m(iv),\quad v\in \mathbb{R}.$$
Hence by the Paley-Wiener theorem, $\Phi$ must be supported in
$[0,\infty)$. A concrete computation (see, page 77 in \cite{fendler})
shows that
$$\int_{0}^{+\infty}e^{-\lambda u}\Phi(u)du=m(\lambda), \quad \lambda\in \mathbb{R}^+.$$
Thus by the spectral theory
\be \int_{0}^{+\infty}e^{-uL}\Phi(u)du &=&
\int_{0}^{+\infty}\Big(\int_{0}^{+\infty}e^{-\lambda
u}dP_\lambda\Big)\Phi(u)du
\\&=&\int_{0}^{+\infty}m(\lambda)dP_\lambda=m(L).
\ee

\noindent By the noncommutative Markov dilation  {Theorem \ref{JXS}},$(T_t)$ admits a Markov dilation. We notice that the dilation property  can  be extended  verbatim to all of $\mathbb{R}$ by using the ultraproduct argument (see \cite[ Corollary 4.3, page 51 and Corollary 4.5, page 54]{jrs} for more details.
Namely, there exists a
larger {finite} von Neumann algebra $\mathcal {N}$, an increasing
filtration $(\mathcal {N}_{s]})_{s\in\mathbb{R}}$ with conditional
expectation $\mathcal {N}_{s]}=\E_{s]}(\mathcal {N})$ and trace
preserving
*-homomorphisms $\pi_s:\M\longrightarrow\mathcal {N}$
such that $\pi_s(\M)\subset\mathcal
{N}_{s]}$ and
\begin{align}
\label{dilation}
\E_{s]}\big(\pi_t(x)\big)=\pi_s(T_{t-s}x),
\quad \quad -\infty< s\le t<\infty,\quad x\in \M.
\end{align}
Where $\pi_{-t}(x):=\beta(\pi_t(x))$ for all $x\in\mathcal {M}$,  $\beta$ is the automorphism of
$\mathcal{N}$ appeared  in  Corollary 4.5 in \cite{jrs}.
Especially, for any $t\in[0,+\infty)$ and $s\in\mathbb{R}$
$$\mathbb{E}_{-s]}\pi_{t-s}x=\pi_{-s}T_tx,\quad \forall x\in \mathcal {M}.$$
{
{
 Consequently, keeping the support of $\Phi$ in $[0,\infty)$  in mind,
for any $x\in \mathcal {M},$
\begin{eqnarray*}
\|m(L)x\|_{L_p(\mathcal M)} &=&
\Big\|\pi_{-s}\int_{0}^{+\infty}\Phi(t)T_txdt\Big\|_{L_p(\mathcal {N})}
=\Big\|\int_{0}^{+\infty}\Phi(t)\mathbb{E}_{-s]}\pi_{t-s}xdt\Big\|_{L_p(\mathcal
{N})}
\\&=&\Big\|\int_{-\infty}^{+\infty}\Phi(t)\mathbb{E}_{-s]}\pi_{t-s}xdt\Big\|_{L_p(\mathcal {N})}=\Big\|\mathbb{E}_{-s]}\int_{-\infty}^{+\infty}\Phi(t)\pi_{t-s}xdt\Big\|_{L_p(\mathcal {N})}
\\&\leq&\Big\|\int_{-\infty}^{+\infty}\Phi(t)\pi_{t-s}xdt\Big\|_{L_p(\mathcal {N})},
\quad \forall s\in \mathbb{R},
\end{eqnarray*}}
which implies that for any $N\in\mathbb{N}$
\begin{eqnarray}\label{m2}
\tau\left(|m(L)x|^p\right)\leq\frac{1}{2N}\int_{-N}^N\tau\Bigg(\Big|\int_{-\infty}^{+\infty}\Phi(t)\pi_{t-s}
xdt\Big|^p\Bigg)ds.
\end{eqnarray}

\noindent For any $M>0$, let $\chi_{(-M-N,\,M+N)}(w)$ be the
characteristic function of $(-M-N,\,M+N)$ and
$f_x(-w)=\pi_wx\chi_{(-M-N,\,M+N)}(w)$, then $f_x(\cdot)\in
L_p(\mathbb{R}, L_p(\mathcal {M}))$.   By our assumption, we obtain
\begin{eqnarray*}
\big\|\Phi*f_x\big\|_{L_p\big(\mathbb{R},L_p(\mathcal {M})\big)}\leq
C\big\|f_x\big\|_{L_p\big(\mathbb{R},L_p(\mathcal {M})\big)}.
\end{eqnarray*}
Enlarging  the integral domain  from $[-N,N]$ to
$\mathbb{R}$ we see that the expression on the right in (\ref{m2}) is
smaller than
\begin{eqnarray*}
&&\frac{1}{2N}\int_{-\infty}^\infty\Bigg(\tau\Big|\int_{-\infty}^{+\infty}\Phi(t)\pi_{t-s}
x\chi_{(-M-N,\,M+N)}(t-s)dt\Big|^p\Bigg)ds\\&=&\frac{1}{2N}\big\|\Phi*f_x\big\|^p_{L_p(\mathbb{R}, L_p(\mathcal
{M}))}\leq\frac{C^p}{2N}\big\|f_x\big\|^p_{L_p(\mathbb{R}, L_p(\mathcal
{M}))}\\&=&\frac{C^p}{2N}\int_{-\infty}^\infty\big\|\pi_{-w}
x\chi_{(-M-N,\,M+N)}(-w)\big\|_{L_p(\mathcal {M})}^pdw
\\&=&\frac{C^p(2N+2M)}{2N}\big\|
x\big\|_{L_p(\mathcal {M})}^p,\quad \quad x\in \mathcal {M}.
\end{eqnarray*}

\noindent Let $N\longrightarrow\infty$, we get
\begin{eqnarray}\label{m5}
\|m(L)x\|_p\leq
C\|x\|_p\quad \quad x\in \mathcal {M}.
\end{eqnarray}

\noindent Since {$\mathcal{M}$} is dense in $L_p(\mathcal {M})$,
$m(L)$ then extends uniquely to a bounded operator
on $L_p(\mathcal {M})$. The
proof is complete.\hfill$\Box$\bigskip

For later use, we record the following corollary on the imaginary
powers of $L$.

\begin{cor}\label{2.3} Suppose that $1<p<\infty$, and that $u\in \R$. Then the operator $L^{iu}$ is bounded
on $L_p(\mathcal {M})$: for any $x$ in $L_p(\mathcal {M})\cap
L_2(\mathcal {M})$, {
\begin{equation}\label{constant}\|L^{iu}x\|_p\leq C
\frac{p^2}{p-1}\Big(1+|u|\Big)\exp(\frac{\pi}{2}|u|)\|x\|_p ,
\end{equation} where $C$ is an absolute constant.}
\end{cor}

\noindent {\bf Proof}  \quad Let $\Phi$ be the distribution with
Fourier transform
$$\hat{\Phi}(v)=(iv)^{iu},\quad v\in \mathbb{R}.$$

\noindent In this case,
$\Phi(t)=\Gamma(-iu)^{-1}t^{-iu-1}$ if $t>0,$ otherwise $\Phi(t)=0,$ where and in the sequel $\Gamma(\cdot)$ is the Gamma function.
$\hat{\Phi}$ is the boundary value of the holomorphic function
$$m(z)=z^{iu}=\exp\big(iu\log|z|-u\arg z\big).$$
We note that $m$ satisfies the {{H\"{o}rmander-Mihlin}} condition:
$$|m(iv)|\le\exp({\pi\over 2}|u|),\ \ \ |{\partial\over\partial v}m(iv)|=|i{u\over v}m(iv)|\le{|u|\over|v|}\exp({\pi\over 2}|u|),\ \ \ v\in\mathbb{R}.$$

\noindent By the H\"{o}rmander multiplier theorem (see \cite{b} or \cite{z}) and the UMD property of noncommutative $L_p(\mathcal
{M})$ space, we deduce that for any $f\in L_p(\mathbb{R},L_p(\mathcal
{M}))$
$$\|\Phi*f\|_{L_p\big(\mathbb{R},L_p(\mathcal {M})\big)}\le C(p)(1+|u|)\exp({\pi\over 2}|u|)\|f\|_{L_p\big(\mathbb{R},L_p(\mathcal {M})\big)}.
$$

\noindent By a result of Parcet {in \cite{parcet}}, we can give an explicit estimate of the constant $C(p)=C[p+{1\over p}]$. {In fact,} let $T$ be the Calder\'{o}n-Zygmund
operator associated to the kernel $\Phi(t)$. Since
$$|\Phi(t)|\le\exp({\pi\over 2}|u|)\lesssim{1\over  |t|},$$
and
$$|\Phi(s)-\Phi(s')|\lesssim{(1+|u|)}\exp({\pi\over 2}|u|)\lesssim\frac{|s-s'|}{|s|^{2}}\ \ \ {\rm{if}}\ \ \ |s-s'|\leq\frac{1}{2}{|s|}.$$
Namely, the kernel $\Phi(t)$ satisfies the size and smoothness conditions with the Lipschitz smoothness parameter $\gamma=1$.
It follows from Theorem A in \cite{parcet} that for any
 $f\in L_p(\mathbb{R},L_p(\mathcal
{M}))$
$$\|\Phi*f\|_{L_p\big(\mathbb{R},L_p(\mathcal {M})\big)}=\|Tf\|_{L_p\big(\mathbb{R},L_p(\mathcal {M})\big)}
\leq C \frac{p^2}{p-1} \Big(1+|u|\Big)\exp(\frac{\pi}{2}|u|)\|f\|_{L_p\big(\mathbb{R},L_p(\mathcal
{M})\big)},$$ where $C$ is a constant independently on $p$. Hence we deduce
from Theorem \ref{main} that \beq\label{iu}\|L^{iu}x\|_p\leq C
\frac{p^2}{p-1}\Big(1+|u|\Big)\exp(\frac{\pi}{2}|u|)\|x\|_p,\ \ \ \ x\in L_p(\mathcal {M})\cap
L_2(\mathcal {M}).\eeq \hfill$\Box$

\begin{remark}\label{good constant} { Note that the operator norm of $L^{iu}$ on
$L_2(\mathcal {M})$ is equal to 1 by the spectral theory. Hence, it is possible to improve the constant in \eqref{constant} by the Riesz-Thorin interpolation theorem. Using verbatim the standard method stated in \cite[Corollary 6.3.1, page 78]{fendler}, we can improve \eqref{constant} as follows which is needed in Section 3
$$\|L^{iu}x\|_p\leq C\frac{p^2}{p-1}\Big(1+{\left|u\right|^{12}}\Big)^{\left|\frac{1}{p}-\frac{1}{2}\right|}\exp\left(\pi\left|\frac{1}{p}-\frac{1}{2}\right||u|\right)\|x\|_p,$$
 where $C$ is an absolute constant. } Moreover, by using Cowling's argument \cite[Corollary 1, page 270]{cowling}, we can  further improve the power index on $|u|$, we leave it to the reader.
\end{remark}

\begin{thm} Suppose that $m$ is a bounded holomorphic
function on the sector $\Sigma_{\phi}$ with $\pi/2<\phi\leq\pi$. Then
for $1<p<\infty,$
$$\|m(L)x\|_p\leq C(p,\phi,m)\|x\|_p,\quad \quad \forall x\in L_p(\mathcal {M}),$$
$C(p,\phi,m)$ is a constant depending only on $p, \phi$ and $m$.
\end{thm}

\noindent {\bf Proof} \quad Let $\Phi$ be the distribution on
$\mathbb{R}$ satisfying
$$\hat{\Phi}(v)=m(iv), \quad \quad v\in \R.$$ Since $m$ extends
analytically to $\Sigma_{\phi}$, $\Phi$ must be supported in
$\R^+\cup \{0\}$. If $v\in\R\setminus\{0\}$, then the disc $D$ with center $iv$ and radius $r$ is contained in
$\Sigma_{\phi}$ provided that
$$r<\sin(\phi-\pi/2)|v|.$$

\noindent By the Cauchy formula, we have
$${dm\over dz}(iv)={1\over 2\pi i}\int_{\partial D}{m(\xi)\over(\xi-iv)^2}d\xi.$$
This implies that
$$|v{\partial \over \partial v}m(iv)|\le {\|m\|_\infty}\text{cosec}(\phi-{\pi\over 2}).$$
Consequently, $\hat{\Phi}(v)$ satisfies the
{{H\"{o}rmander-Mihlin}} condition. Since the noncommutative
$L_p(\mathcal {M})$ space has UMD property, we claim that for
$1<p<\infty,$
$$\|\Phi*f\|_{L_p\big(\mathbb{R},L_p(\mathcal {M})\big)}\leq C(p,\phi,m)\|f\|_{L_p\big(\mathbb{R},L_p(\mathcal {M})\big)},\ \ \ f\in{L_p\big(\mathbb{R},L_p(\mathcal {M})}.$$
The desired result immediately follows from Theorem \ref{main}.\hfill$\Box$\bigskip

The theorem above shows that $L$ admits a bounded
$H^\infty(\Sigma_\phi)$ functional calculus with $\pi/2<\phi\leq\pi$.
By a standard angle reduction principle for noncommutative
semigroup, see \cite[Proposition 5.8]{JLX}, $L$ actually admits a
bounded $H^\infty(\Sigma_\phi)$ functional calculus for any
$\phi>|\frac{1}{p}-\frac{1}{2}|\pi$, which positively
answers the question raised in \cite[Remark 5.9]{JLX}. Then we summarize the main result
of this section as follows.

\begin{thm}\label{th2.4}  Suppose that $m$ is a bounded holomorphic function on the sector
$\Sigma_\psi$. If $|\frac{1}{p}-\frac{1}{2}|\pi<\psi,$ then for
$1<p<\infty$, \be\|m(L)x\|_p\leq C(p,\psi,m)\|x\|_p,\quad \quad x\in
L_p(\mathcal {M}),\ee
for some constant $C(p,\psi,m)$.
\end{thm}

\begin{rk}
By tensoring $\mathcal {M}$ with $M_n$,   the
algebra of $n\times n$ matrices, for any $n$, the theorem above implies that $L$ has a completely bounded
$H^\infty$ functional calculus in $\Sigma_\psi$ with $\phi>|\frac{1}{p}-\frac{1}{2}|\pi$. We refer the interested
reader to \cite{JLX} for more information on $H^\infty$ functional calculus.
\end{rk}

Since every bounded $H^\infty$ functional calculus implies square function estimates, we have the following corollary from Theorem \ref{th2.4}. We refer to
\cite[Theorem 7.6 or Corollaries 7.7 and 7.10]{JLX} for more details on square functions.

\begin{cor}\label{2.6}
Suppose that $m$ is a bounded holomorphic function on the sector $\Sigma_\psi$ with $\psi>|\frac{1}{p}-\frac{1}{2}|\pi$.

\noindent Then for any $x\in L_p(\mathcal{M})$ the following hold:
\begin{itemize}
\item[{\em (1)}] For $1<p<2$,
\begin{align*}
\|x\|_{ L_p(\mathcal{M})}\thickapprox\inf
&\left\{\ \
\left\|\left(\int_0^\infty t\left|{\partial\over\partial t}(T_t(x_1))\right|^2dt\right)^{1/2}\right\|_{ L_p(\mathcal{M})}+\right.\\
&\left.\ \ \ \ \ \left\|\left(\int_0^\infty t\left|{\partial\over\partial t}(T_t(x_2))^*\right|^2dt\right)^{1/2}\right\|_{ L_p(\mathcal{M})}\ \ \right\},
\end{align*}
where the infimum runs over all $x_1, x_2\in L_p(\mathcal{M})$ such that $x=x_1+x_2$.

\item[{\em (2)}] For $2\le p<\infty$,
\begin{align*}
\|x\|_{L_p(\mathcal{M})}\thickapprox\max
&\left\{\ \
\left\|\left(\int_0^\infty t\left|{\partial\over\partial t}(T_t(x))\right|^2dt\right)^{1/2}\right\|_{L_p(\mathcal{M})}, \right.\\
&\left.\ \ \ \ \ \left\|\left(\int_0^\infty t\left|{\partial\over\partial t}(T_t(x))^*\right|^2dt\right)^{1/2}\right\|_{ L_p(\mathcal{M})}\right\}.
\end{align*}
\end{itemize}\end{cor}

\begin{rk}
{It is worth noticing that a similar result is obtained in \cite[pp. 68]{JLX} with a different method. We should emphasize here that the ours
 is much simpler and works without the hypothesis of $H^\infty$-functional
calculus for $L$. }

\end{rk}

%
\section{Noncommutative maximal inequalities}

We now turn to maximal inequalities.
We first recall the definition of noncommutative maximal functions
introduced by Pisier  (see \cite{pisier}) and Junge (see
\cite{junge}). Let $1\le p\leq\infty$,
we define $L_p(\mathcal {M},\ell_\infty)$ to be the space of all sequences $x=(x_n)_{n\geq1}$
in $L_p(\mathcal {M})$, which admit a factorization of the following form:
there exist $a,b\in L_{2p}( \mathcal {M})$  and a bounded sequence $y=(y_n)\subset L_\infty(\mathcal {M})$ such that
$$x_n=ay_nb,\ \ \ n\ge 1.$$
The norm of $x$ in $L_p(\mathcal {M},\ell_\infty)$ is given by
$$
\|x\|_{L_p(\mathcal {M},\ell_\infty)}=\inf\left\{\big\|a\big\|_{{2p}}\sup_n\big\|y_n\big\|_\infty\big\|b\big\|_{{2p}}\right\},$$
where the infimum is taken over all factorizations of $x$ as above.
It is easy to see that $L_p(\mathcal {M},\ell_\infty)$ is a Banach space with the norm $\|\cdot\|_{L_p(\mathcal {M},\ell_\infty)}$,
and a positive sequence
$x=(x_n)$ belongs to $L_p(\mathcal {M},\ell_\infty)$ if and only if there is $a\in L_p^+(\mathcal {M})$
such that $x_n\leq a$ for all $n$. Moreover, in this case,
$$\|x\|_{L_p(\mathcal {M}, \ell_\infty)}=\inf\Big\{\|a\|_p: a\in L_p^+(\mathcal {M})\ \text{such\ that}\ x_n\leq a
,\,\forall n\geq1\Big\}.$$

\noindent The norm of $x$ in $L_p(\mathcal {M},\ell_\infty)$ is conventionally denoted by $\|\sup_{n\geq1}^+x_n\|_p$. Please
note that $\|\sup_{n\geq1}^+x_n\|_p$ is
just a notation since $\sup_{n\geq1}x_n$ does not make any
sense in the noncommutative setting. We use this notation only for convenience.

\begin{rk} The definition of $L_p(\mathcal {M}, \ell_\infty)$ can be
extended to an arbitrary index set  $I$.  Then
$L_p(\mathcal {M},\ell_\infty(I))$  can be defined similarly as
before. More precisely, $L_p(\mathcal {M},\ell_\infty(I))$ consists of
all families $(x_i)_{i\in I}$ in $L_p(\mathcal {M})$ which can be factorized as $x_i=ay_ib$
with $a,b\in L_{2p}(\mathcal {M})$ and a bounded family $(y_i)_{i\in I}\subset L_\infty(\mathcal {M})$.
The norm of $(x_i)_{i\in I}$ in $L_p(\mathcal {M}, \ell_\infty(I))$ is
defined as
$$\inf\left\{\big\|a\big\|_{2p}\sup _i\big\|y_i\big\|_\infty\big\|b\big\|_{2p}\right\},$$
the infimum running over all factorizations as above. As before, this norm is also
denoted by $\|\sup_{i\in I }^+x_i\|_{p}.$
\end{rk}

\begin{rk} One can easily check that for any index set $I$ and $1\leq p\leq \infty$,  a family $(x_i)_{i\in I}$
in $L_p(\mathcal {M})$ belongs to $L_p(\mathcal {M},\ell_\infty(I))$ if and only if $$\sup_{J\subset I, J \ \text{is\ a\ finite\ set}}\ \|\sup_{i\in J }\!^+x_i\|_{p}<\infty.$$ If this is the case, then
$$\|\sup_{i\in I }\!^+x_i\|_{p}=\sup_{J\subset I, J \ \text{is\ a\ finite\ set}}\ \|\sup_{i\in J }\!^+x_i\|_{p}.$$
\end{rk}

The main result of this section is
relevant to the noncommutative maximal function $\|\sup_{z\in
\Sigma_\psi}^+T_zx\|,$ where $\Sigma_\psi$ is a sector in $\mathbb{C}$, which
generalizes the Theorem 5.1 and
Corollary 5.11 in \cite{JX}, and Corollary 5.7  in \cite{jm1}. Moreover,
we will see in the next section that our maximal inequality implies that $(T_zx)_z$
converges bilaterally almost uniformly for any $x\in L_p(\mathcal {M})$. In addition, for $p>2$ the bilateral almost uniform
convergence can be improved to the almost uniform convergence.
For {formulating} this result we need  further notation from \cite{DJ}. Let $2\leq p\leq\infty$ and $I$ be an index set.
We define the space $L_p(\mathcal {M},\ell^c_\infty(I))$ as  the family of all $(x_i)_{i\in I}\subset L_p(\mathcal {M})$
for which there are an $a\in L_p(\mathcal {M})$ and $(y_i)_{i\in I}\subset L_\infty(\mathcal {M})$ such that
$$x_i=y_ia\ \ \ \ \text{and}\ \ \ \ \sup_{i\in I}\|y_i\|_\infty<\infty.$$
$\|(x_i)\|_{L_p(\mathcal {M},\ell^c_\infty(I))}$ is then defined to be the infimum $\{\sup_{i\in I}\|y_i\|_\infty\|a\|_p\}$
over all factorizations of $(x_i)$ as above. It is easy to check that $\|\cdot\|_{L_p(\mathcal {M},\ell^c_\infty(I))}$
 is a norm, which makes $L_p(\mathcal {M},\ell^c_\infty(I))$ a Banach space. Note that $(x_i)\in L_p(\mathcal {M},\ell^c_\infty(I))$ if and only if $(x^*_ix_i)\in L_{p\over 2}(\mathcal {M},\ell_\infty(I))$. If $I=\mathbb{N} $, $L_p(\mathcal {M},\ell^c_\infty(I))$ is simply denoted by
 $L_p(\mathcal {M},\ell^c_\infty)$.
 To state our maximal inequalities we also need the following lemma.


\begin{lem}{\cite{cowling}}\label{lemma2}  Let \beq\label{m}
m_\theta(\lambda)=\exp(-e^{i\theta}\lambda)-\int_0^1\exp(-t\lambda)dt,\eeq
where $|\theta|\leq\pi/2,$ and $n_\theta=m_\theta\circ exp$. \ Then for any $u\in\mathbb{R}$
$$\hat{n}_\theta(u)=\Big(e^{-\theta u}-(1+iu)^{-1}\Big)\Gamma(-iu) \ \ \ \  \ \text{and}\ \ \ \ \ \ |\hat{n}_\theta(u)|\lesssim\exp\Big(\left[|\theta|-\frac{\pi}{2}\right]|u|\Big).$$
 \end{lem}

\begin{thm} \label {max}  Suppose that $1<p<\infty$, and $$0\leq\psi/\pi<1/2-\left|1/p-1/2\right|.$$
Let ${\Sigma}_\psi$ be the sector $\{z\in \mathbb{C}:|\arg z|<
\psi\}$. Then there exists a constant $C$ depending only on $p$ and $\psi$ such that
\begin{eqnarray}
\label{th3.3}
\|\sup_{z\in{\Sigma}_\psi}\!\!^+T_zx\|_p\leq C\big\|x\big\|_p,\quad \quad x\in L_p(\mathcal {M}).
\end{eqnarray}
Moreover, if $p>2$, then
\begin{eqnarray}
\label{cor3.4}
\|(T_zx)_{z\in\Sigma_\psi} \|_{L_p(\mathcal {M},\ell^c_\infty({\Sigma}_\psi))}\leq C  \|x \|_p,\ \ \ \ \forall x\in L_p(\mathcal {M}).
\end{eqnarray}
 \end{thm}

\noindent{\bf Proof } \ Let $n_{\theta}(t)=m_{\theta}\big(\exp(t)\big), t\in\mathbb{R}$ and
\beq\label{3.2}
m_{\theta}(\lambda)=\exp(-e^{i\theta}\lambda)-\int_0^1\exp(-t\lambda)dt,\quad
\lambda\in\R^+,\eeq
where $|\theta|<\psi.$
It follows from Fourier transform that
$$m_{\theta}(\lambda)=\frac{1}{2\pi}\int_{-\infty}^\infty\hat{n}_{\theta}(u)\lambda^{iu}
du,\quad \quad \lambda\in\mathbb{R}^+.$$

\noindent  By functional
calculus we have,
\beq\label{3.3}
m_{\theta}(tL )x=\frac{1}{2\pi}\int_{-\infty}^\infty\hat{n}_{\theta}(u)t^{iu}L^{iu}x du,\ \ \ \ \ t\in \mathbb{R}^+.
\eeq

\noindent Let $z=te^{i\theta}$ with $|\theta|<\psi$ and $t>0$, then it follows from
\eqref{3.2}, by functional calculus,  that
\begin{eqnarray}
\label{3.4}
T_zx=e^{-zL}x&=&m_{\theta}(tL)x+\int_0^1\exp(-st
L)xds\notag\\
&=&m_{\theta}(tL)x+\frac{1}{t}\int_0^tT_sxds.
\end{eqnarray}

\noindent Consequently,
\beq\label{a}\displaystyle\big\|\sup_{z\in{\Sigma}_\psi}\!\!^+T_zx\big\|_p
\leq\big\|\sup_{t>0,|\theta|<\psi}\!\!\!\!\!\!^+ m_{\theta}(tL)x\big\|_p
+\big\|\sup_{t>0}\!^+\frac{1}{t}\int_0^tT_sxds\big\|_p.\eeq

\noindent Since
\begin{eqnarray}
\label{mtheta}
\big\|\displaystyle\sup_{t>0,|\theta|<\psi}\!\!\!\!\!\!\!^+ m_{\theta}(tL)x\big\|_p
&=&\big\|\displaystyle\sup_{t>0,|\theta|<\psi}\!\!\!\!\!\!\!^+ \ \ \displaystyle\frac{1}{2\pi}\int_{-\infty}^\infty\hat{n}_{\theta}(u)t^{iu}L^{iu}x du\big\|_p\nonumber\\
&\le&\displaystyle\frac{1}{2\pi}\int_{-\infty}^\infty\big\|\displaystyle\sup_{t>0,|\theta|<\psi}\!\!\!\!\!\!\!^+ \ \hat{n}_{\theta}(u)t^{iu}L^{iu}x \big\|_pdu\nonumber\\
&\le&\displaystyle\frac{1}{2\pi}\int_{-\infty}^\infty\displaystyle\ \ \sup_{t>0,|\theta|<\psi}\bigg|
{\hat{n}_{\theta}(u)t^{iu}\over\exp([\psi-{\pi\over 2}]|u|)}\bigg|\ \bigg\|\exp\big([\psi-{\pi\over 2}]|u|\big)L^{iu}x \bigg\|_pdu\\
&\lesssim&\Big(\frac{1}{2\pi}\int_{-\infty}^\infty\exp\big([\psi-{\pi\over 2}]|u|\big)\parallel\mid
L^{iu}\parallel\mid\textcolor{red}{_p} du\Big)\big\|x\big\|_p\nonumber\\
&=&C(p,\psi)\big\|x\big\|_p,\nonumber
\end{eqnarray}

\noindent where $\textcolor{red}{\parallel\mid L^{iu}\parallel\mid_p}$ is the operator norm of
$L^{iu}$ on $L_p(\mathcal {M})$ and $$C(p,\psi)=\frac{1}{2\pi}\int_{-\infty}^\infty\exp([\psi-{\pi\over 2}]|u|)\parallel\mid
L^{iu}\parallel\mid\textcolor{red}{_p} du.$$
{ It follows from Corollary \ref{2.3} and Remark \ref{good constant}, we have
\begin{eqnarray*}
C(p,\psi)&\lesssim&
\frac{p^2}{p-1}\frac{1}{2\pi}\int_{-\infty}^\infty
\Big(1+{\left|u\right|^{12}}\Big)^{\left|\frac{1}{p}-\frac{1}{2}\right|}\exp\left(\pi\left|\frac{1}{p}-\frac{1}{2}\right||u|\right)\exp\Big([\psi-\frac{\pi}{2}]|u|\Big)du <\infty,\end{eqnarray*} where the finiteness of last integral can be found in \cite[page 81]{fendler}.
Thus we deduce that}
\beq \label{b}
\big\|\sup_{t>0,|\theta|<\psi}\!\!\!\!\!^+ \ m_{\theta}(tL)x\big\|_p\lesssim
C(p,\psi)\|x\|_p, \ \ \  \forall x\in L_p(\mathcal {M}).
\eeq

\noindent Similarly, we also get that
\beq \label{bb}
\bigg\| (m_{\theta}(tL)x)_{t>0,|\theta|<\psi}\bigg\|_{L_p(\mathcal {M},\ell^c_\infty(\Sigma_\psi)))}\lesssim C(p,\psi)\|x\|_p, \ \ \ \forall x\in L_p(\mathcal {M}).\eeq On the other hand,
Theorem 4.5 in \cite{JX} implies that
\beq\label{c}
\bigg\|\sup_{t>0}\!^+\frac{1}{t}\int_0^tT_sxds\bigg\|_p\leq C_p\|x\|_p,\ \ \forall x\in L_p(\mathcal {M}) \ \ \ \text{for}\ \ \ 1<p<\infty,
\eeq

\noindent and
\beq\label{e}
\bigg\| \left(\frac{1}{t}\int_0^tT_sxds\right)_{t>0}\bigg\|_{L_p(\mathcal {M},\ell^c_\infty(\mathbb{R}^+))}\leq C_p\|x\|_p,\ \ \forall x\in L_p(\mathcal {M}) \ \ \ \text{for}\ \ \ 2<p<\infty.
\eeq

\noindent Combining \eqref{a}, \eqref{b} and \eqref{c} implies the desired estimate (\ref{th3.3}).
And the desired estimate (\ref{cor3.4}) follows from \eqref{a}, \eqref{bb} and \eqref{e}.  \hfill$\Box$\bigskip

\section{Individual Ergodic Theorems}

In this section{, motivated by  {Proposition 7}  in \cite{xu}, }
we apply the maximal inequalities proved in the previous section to the pointwise ergodic convergence.
To this end we need an appropriate analogue for the noncommutative setting of the usual almost everywhere convergence. This is the almost
uniform convergence introduced by Lance in \cite{L}.

\begin{defi} Let $\mathcal {M}$ be a von Neumann algebra equipped
with a {finite} normal faithful trace $\tau.$ Let $x_n,x\in
L_0(\M)$.

\begin{itemize}
\item[\em(1)] \ $(x_n)$ is said to converge bilaterally almost uniformly
(b.a.u. in short) to $x$ if for every $\varepsilon>0$ there is a
projection $e\in \M$ such that
$$\tau(e^\bot)<\varepsilon \quad \quad {\rm and }\quad \quad
\lim _{n\rightarrow\infty}\|e(x_n-x)e\|_\infty=0.$$

\item[\em(2)] \ $(x_n)$ is said to converge almost uniformly (a.u. in short) to
$x$ if for every $\varepsilon>0$ there is a projection $e\in \M$
such that
$$\tau(e^\bot)<\varepsilon \quad \quad {\rm and }\quad \quad
\lim _{n\rightarrow\infty}\|(x_n-x)e\|_\infty=0.$$
\end{itemize}\end{defi}

In the commutative case, both convergences in the definition above
are equivalent to the usual almost everywhere convergence by virtue
of Egorov's theorem.  However they are different in the
noncommutative case. Similarly, we {can} introduce these notions of convergence for functions with values in
$L_0(\mathcal {M})$ and for nets in $L_0(\mathcal {M})$.

\begin{thm} \label {4.2} Suppose that $1<p<\infty$, and
that$$0\leq\psi/\pi<1/2-|1/p-1/2|.$$
Then for any   $x\in  L_p(\mathcal {M})$ the following hold:

\begin{itemize}
\item[{\em (1)}] The operators $T_zx$ converge bilaterally almost uniformly to
$x$ for $1<p\leq 2$ and almost uniformly to $x$ for $2< p<\infty$ as
$z$ tends to 0 in ${\Sigma}_\psi.$

\item[{\em (2)}] The operators $T_zx$ converge bilaterally almost uniformly to
$P(x)$ for $1<p\leq 2$ and almost uniformly to $P(x)$ for $2< p<\infty$ as
$z$ tends to $\infty$ in ${\Sigma}_\psi,$ where $P$ denotes the projection from
$  L_p(\mathcal {M})$ onto the fixed point subspace of the semigroup $(T_t)_{t>0}$.
\end{itemize}\end{thm}

\noindent{\bf Proof}\quad (1) \  Let   $x\in L_2(\mathcal {M})\cap
L_p(\mathcal {M})$ and $s>0.$ Let $D$ be the disc of center $s$ and radius
$r$ with $r<s\sin\psi$. For  any $z\in D$, by the
vector-valued Cauchy formula, we have
$$T_z(x)=\frac{1}{2\pi i}\int_{\partial D}\frac{T_\zeta(x)d\zeta}{\zeta-z}.$$

\noindent Thus$$T_z(x)-T_s(x)=\frac{z-s}{2\pi
i}\int_{\partial D}\frac{T_\zeta(x)d\zeta}{(\zeta-z)(\zeta-s)}.$$

\noindent By the convexity of the operator valued function: $x\longmapsto |x|^2$,
\begin{eqnarray}\label{4.0}|T_z(x)-T_s(x)|^2\leq\frac{|z-s|^2}{4\pi^2
}\int_{\partial D}\frac{|d\zeta|}{|(\zeta-z)(\zeta-s)|^2}
\int_{\partial D}|T_\zeta(x)|^2|d\zeta|\leq C|z-s|^2a
\end{eqnarray}

  \noindent for $|z-s|<{r\over 2}, $ where $C$ denotes
a positive constant independent of $z$ and
$a=\int_{\partial D}|T_\zeta(x)|^2|d\zeta|.$
Note that $a\in L_1(\mathcal {M})$. It follows that there exists a
contraction $u\in L_\infty(\M)$ (depend on $z$ and $s$) such that
$$T_z(x)-T_s(x)=C|z-s|ua^{1/2}.$$
For any  $\varepsilon>0$, let $e=e_{(0,{\|a^{1/2}\|_2/\epsilon^{1/2}})}(a^{1/2})$ be the spectral projection of
$a^{1/2}$ on the interval $(0,{\|a^{1/2}\|_2/\varepsilon^{1/2}})$. Then
$$\tau(e^\bot)<\varepsilon \quad \quad {\rm and }\quad \quad
a^{1/2}e\in L_\infty(\M).$$
Therefore $$\Big\|\Big(T_z(x)-T_s(x)\Big)e\Big\|_\infty=C|z-s|\big\|ua^{1/2}e\big\|_\infty\leq
C|z-s|\big\|a^{1/2}e\big\|_\infty.$$ Consequently,
$$\lim_{z\rightarrow s}\Big\|\Big(T_z(x)-T_s(x)\Big)e\Big\|_\infty=0.$$
Namely, $\lim_{z\rightarrow s}T_z(x)=T_s(x)$ almost uniformly. It
then follows that $\lim_{z\rightarrow 0}T_z\big(T_s(x)\big)=T_s(x)$
almost uniformly for all $x\in L_2(\mathcal {M})\cap L_p(\mathcal
{M})$. Since the linear span of $ \big\{T_sx:x\in L_2(\mathcal {M})\cap
L_p(\mathcal {M}),s>0\big\}$ is dense in $L_p(\mathcal {M})$,   our desired results then follows from
Theorem \ref{max}.

Indeed, for $x\in L_p(\mathcal {M})$ and $\epsilon>0$, take a sequence $(x_n)$ in the span of $ \big\{T_sx:x\in L_2(\mathcal {M})\cap
L_p(\mathcal {M}),s>0 \big\}$ such that
$$\|x_n-x\|_p<{1\over 2^n}\left({\epsilon\over 2^n}\right)^{1\over p}.$$

\noindent Let $e_{1,n}=e_{(0,1/2^n)}(|x-x_n|)$ be the spectral
projection of $|x-x_n|$ on the interval $(0,1/2^n)$, then
$$\tau(1-e_{1,n})\leq\frac{\|x-x_n\|^p_p}{(1/2^n)^p}<\frac{\varepsilon}{2^n}\quad {\rm and}\quad
\| |x-x_n|e_{1,n}\|_\infty<\frac{1}{2^n}.$$
Set $e_1=\bigwedge_n e_{1,n}$. We have
\begin{eqnarray}\label{4.1}
\tau(1-e_1)<\varepsilon \quad \quad {\rm and }\quad \quad
\| |x-x_n|e_{1}\|_\infty<\frac{1}{2^n}.
\end{eqnarray}

\noindent From inequality (\ref{th3.3}) in Theorem \ref{max}, we know that
$$\big\|\sup_{z\in{\Sigma}_\psi}\!\!^+T_z(x-x_n)\big\|_p\leq C(p,\psi)\big\|x-x_n\big\|_p.$$
That is, there are $a,b\in L_{2p}(\mathcal{M})$ and $(c_z)\subset L_\infty(\mathcal{M})$ such that
$$\sup_{z\in\Sigma_\psi}\|c_z\|_\infty\le 1,\ \ \|a\|_{2p}=\|b\|_{2p}\lesssim\|x-x_n\|_p^{1\over 2},$$
and
$$T_z(x-x_n)=ac_zb.$$

\noindent Let $e_{2,n}=e_{(0,(1/2)^{n\over 2})}(|a|)\bigwedge e_{(0,(1/2)^{n\over 2})}(|b|)$,
then
$$\tau(1-e_{2,n})\leq\frac{\|a\|^{2p}_{2p}}{(1/2^n)^{p}}\lesssim\frac{\varepsilon}{2^n}\quad {\rm and}\quad
\|e_{2,n}T_z(x-x_n)e_{2,n}\|_\infty\leq\frac{1}{2^n}.$$

\noindent Set $e_2=\bigwedge_n e_{2,n}$. We have
\begin{eqnarray}\label{4.2}
\tau(1-e_2)\lesssim\varepsilon \quad \quad {\rm and }\quad \quad
\| e_{2}T_z(x-x_n)e_{2}\|_\infty\leq\frac{1}{2^n}.
\end{eqnarray}

\noindent Since $\lim_{z\rightarrow 0}T_z\big(x_n\big)=x_n$
almost uniformly for all $n,$
there is a projection $e_{3,n}$ such that
$$\tau(1-e_{3,n})<\frac{\varepsilon}{2^n}\ \ \ \ \text{and}\ \ \  \ \lim_{\Sigma_\psi\ni z\to 0 }\big\| (T_zx_n-x_n)e_{3,n}\big\|_\infty=0..$$ Let $e_3=\bigwedge_n
e_{3,n}$, then
\begin{eqnarray}\label{4.3}\tau{(1-e_3)}<\varepsilon\quad \quad {\rm and }
\quad
\quad\lim_{z\in\Sigma_\psi }\big\| (T_zx_n-x_n)e_3\big\|_\infty=0.
\end{eqnarray}

\noindent Take
$e=e_1\wedge e_2\wedge e_3,$
then
\begin{eqnarray*} \tau(e^\bot)\lesssim \varepsilon,
\end{eqnarray*}
 and
\begin{eqnarray*}
\|e(T_zx-x)e\|_\infty&\leq&\|e(T_zx-T_zx_n)e\|_\infty+\|e(T_zx_n-x_n)e\|_\infty+\|e(x_n-x)e\|_\infty.
\end{eqnarray*}
Thus it follows from formulas (\ref{4.1})-(\ref{4.3}) that
$$\lim_{\Sigma_\psi\ni z\to 0}\|e(T_zx-x)e\|_\infty=0.$$
Thus the first part of (1) is proved. The case for $p>2$ can be similarly proved by using inequality (\ref{cor3.4}). This completes the proof of (1).

(2) \ Now we turn to the second part of the theorem. First $L$ induces a canonical splitting on $ L_p(\mathcal {M})$ for $1<p<\infty:$ $L_p(\mathcal {M})=N(L)\oplus \overline{R(L)}$. Moreover, $N(L)$ is the fixed point subspace of the semigroup $(T_t)_{t>0}$. Thus it suffices to prove that for any $x\in \overline{R(L)}$,
 $T_zx$ converge bilaterally almost uniformly to $P(x)$ for $1<p\leq 2$ when $z\rightarrow \infty$. Using (3.2) in Theorem 3.4 as in the previous part of the proof, we need only to do
 this for $x$ in a dense subset of $\overline{R(L)}.$ It is well known that $\{T_{t+s}(y)-T_{s}(y):s>0,t>0, y\in L_p(\mathcal {M})\}$ is such a subset.
 Thus we are reduced to prove the above convergence for $x=T_{t+s}(y)-T_{s}(y).$ Let ${\tilde{\Gamma}}$ be the boundary of ${\Sigma}_\psi.$ Noting that $\|\lambda(\lambda-L)^{-1}\|$is bounded on ${\tilde{\Gamma}}$, we have the integral  representation of $T_z$
 $$T_z=\frac{1}{2\pi i}\int_{{\tilde{\Gamma}}}e^{-z\lambda}(\lambda-L)^{-1}d\lambda,\quad z\in \Sigma_\psi.$$ Then for $x=T_{t+s}(y)-T_{s}(y),$
 $$T_z(x)=\frac{1}{2\pi i}\int_{{\tilde{\Gamma}}}\Big(e^{-(z+t+s)\lambda}-e^{-(z+s)\lambda}\Big)(\lambda-L)^{-1}(y)d\lambda.$$
 Again by the convexity of the operator valued function: $x\longmapsto |x|^2$,
\begin{eqnarray}
\label{4.5}
|T_z(x)|^2&\leq&\frac{t^2}{4\pi^2}\int_{{\tilde{\Gamma}}}|e^{-2z\lambda}||d\lambda|\int_{{\tilde{\Gamma}}}|e^{-2(t+s)\lambda}||(\lambda-L)^{-1}(y)|^2|d\lambda|
\\&=& C\frac{1}{|z|}a,\nonumber
\end{eqnarray}
where $C$ denotes a positive constant independent of $z$ and $a$ is a positive operator in $L_1(\mathcal {M}).$
The  {remaining part of the proof}
  is similar to that of (1)
  and we omit it here.
\hfill$\Box$


\bigskip

\noindent {\bf Acknowledgements} \
{The authors thank the anonymous reviewer for useful suggestions and comments which improve the final version. The authors also wish to acknowledge  Professor Quanhua Xu for his invaluable guidance. The
Part of this research was performed during the second  author's visit to
Laboratoire de Math\'{e}matiques, Universit\'{e} de Franche-Comt\'{e}. He thanks the institution
for their hospitality and support.}

%
%


\bigskip

\begin{thebibliography}{999}

\bibitem{b} Burkholder D.,
\newblock Martingales and singular integral in Banach spaces,
\newblock {\em  Handbook of Geometry of Banach spaces}, Vol.1: 233-269, North-Holland, 2001.

\bibitem{CLS} Chilin V., Litvinov S., Skalski A., A few remarks in non-commutative ergodic theory, J. Operator Theory, 2005, 53: 331-350.

\bibitem{crw} Coifman R., Rochberg R., Weiss G.,
\newblock Applications of trnasference: the $L^p$-version von Neuman's inequality and the Littlewood-Paly-Stein theory,
\newblock {\em  Linear Spaces and Approximation}, 53--67, Basel, 1978.

\bibitem{cw} Coifman R., Weiss G.,
\newblock Transference methods in analysis,
\newblock {\em  CBMS regional conference series in mathematics, No.31, A.M.S., Providence, R.I.}, 1976.

\bibitem{cowling} Cowling M.,
Harmonic analysis on semigroups,
Ann. Math., 1983, 117: 267--283.

\bibitem{cl} Cowling M., Leinert M.,
Pointwise convergence and semigroups acting on vector-valued functions,
Bull. Aust. Math. Soc., 2011, 84: 44-48.

\bibitem{D}
Dabrowski Y., A non-commutative path space approach to stationary
free stochastic differential equations, arXiv:1006.4351.

\bibitem{Dab}
{
Dabrowski Y., A free stochastic partial differential equation,
Ann. Inst. H. Poincar¨¦ Probab. Statist., 2014, 50: 1404-1455.}


\bibitem{DJ} Defant A., Junge M.,  Maximal theorems of Menchoff-Rademacher type in noncommutative $L_q$ spaces,
J. Funct. Anal., 2004, 206: 322-355.

\bibitem{fendler} Fendler G.,
On dilation and transference for continuos one-parameter semegroup of
positive contraction on $L^p$-spaces,
Ann. Univ. Sarav. Ser. Math., 1998, 9(1), iv+97 pp.

\bibitem{GL} Goldstein M., Litvinov S., Banach principle in the space of $\tau$-measurable operators,
Studia Math., 2000, 143: 33-41.

\bibitem{hansen} Hansen F., An operator inequality,
Math. Ann.,  1979/80, 246: 249-250.

\bibitem{h} Hyt\"{o}nen T.,
Littlewood-Paley-Stein theory for semigroups in UMD
spaces, Rev. Mat. Iberoam., 2007, 23: 973-1009.

\bibitem{junge} Junge M.,
\newblock Doob's inequality for noncommutative martingale,
J.Reine Angew. Math., 2002, 549: 149-190.

\bibitem{JLX} Junge, M., Le Merdy, C., Xu, Q.
\newblock $H^\infty$ functional calculus and square functions on noncommutative $L_p$-spaces,
\newblock {\em  Ast¨¦risque 305}, 2006.

\bibitem{jm2}
Junge M., Mei T.,
\newblock Noncommutative Riesz transforms-a probabilstic approach,
Amer. J. Math., 2010, 132: 611-681.

\bibitem{jm1} Junge M., Mei T.,
\newblock BMO spaces associated with semigroups of operators,
Math. Ann., 2012, 352: 691-743.

\bibitem{jrs} {Junge M., Ricard E., Shlyakhtenko D.,
\newblock {Noncommutative diffusion semigroups and free probability, }
 {  to appear.}}

\bibitem{JX} Junge M.,  Xu Q.,
Noncommutative maximal ergodic theorems,
J. Amer. Math. Soc., 2007, 20(2): 385-439.

\bibitem{K} Kriegler C.
Analyticity angle for non-commutative diffusion
semigroups, J. London Math. Soc., 2011,   83: 168-186.

\bibitem{ks} Kunstmann P., \v{S}trkalj, \v{Z}:
\newblock $H^\infty$-calculus for submarkovian generators, Proc. Amer. Math.
Soc., 2003, 131(7):   2081-2088.

\bibitem{L} Lance E.,
Ergodic theorems for convex sets and operator algbras,
Invent.Math., 1976, 37(3): 201-214.

\bibitem{mtx} Martineza T., Torrea J., Xu Q.,
\newblock Vector-valued Littlewood-Paley-Stein theory for semigroups, Adv.
Math., 2006, 203:  430-475.

\bibitem{m} Mei T.,
\newblock Tent spaces associated with semigroup of operators,
 J. Funct. Anal., 2008, 255: 3356-3406.



\bibitem{mp}Mei T., Parcet J.,
\newblock Pseudo-localization of singular integral and noncommutative Littlewood-Paley inequalities,
 Int. Math. Res. Not., 2009, 8: 1433-1487.

\bibitem{parcet}Parcet J.,
\newblock Pseudo-localization of singular integral and noncommutative Calder$\acute{o}$n-Zygmund theory,
J. Funct. Anal., 2009, 256: 509-593.

\bibitem{pisier} Pisier, G.,
\newblock Noncommutative vector-valued $L^p$-spaces and completely $p$-summing maps,
\newblock {C. R. Acad. Sci. Paris S¨¦r. I Math., 1993, 316(10): 1055-1060. }

\bibitem{pxhand} Pisier,G., Xu Q.,
\newblock Noncommutative $L^p$ spaces,
\newblock {\em Handbook of Geometry of Banach spaces,} 2003, 2: 1459-1517.

\bibitem{Stein} Stein, E.,
\newblock Topics in Harmonic analysis related to the Littlewood-Paley Theory,
\newblock {\em  Annals of Mathematics Studies, No.63, Princeton
University Press}, 1970.

\bibitem{taggart} Taggart R., Pointwise convergence for semigroups in vector-valued
$L_p$ spaces, Math. Z., 2009, 261: 933-949.

\bibitem{xu} {Xu Q., \ $H^\infty$ functional calculus and maximal inequalities for
semigroups of contractions on vector-valued  $L_p$-spaces,
\newblock { Int. Math. Res. Not., 2014, rnu104, 18 pages, doi:10.1093/imrn/rnu104.}}


\bibitem{z} Zimmermann F., On vector-valued Fourier multiplier theorems, Studia
Math., 1989, 93: 201-222.


\end{thebibliography}
\end{document}